\documentclass[11pt]{article}
\usepackage[T1]{fontenc}
\usepackage{lmodern}
\usepackage[a4paper,hmargin=28mm,vmargin=26mm]{geometry}
\usepackage{amsmath,amssymb,amsthm}
\usepackage{microtype}
\usepackage[hidelinks]{hyperref}

\newtheorem{theorem}{Theorem}
\newtheorem{proposition}[theorem]{Proposition}
\newtheorem{corollary}[theorem]{Corollary}
\theoremstyle{definition}
\newtheorem{example}[theorem]{Example}

\newcommand{\Q}{\mathcal Q}
\newcommand{\Z}{\mathbb Z}
\newcommand{\supp}{\operatorname{supp}}
\newcommand{\bt}{\boldsymbol t}
\newcommand{\ind}{\mathbf 1}
\newcommand{\conv}{\operatorname{conv}}

\title{Gamma-positivity of support enumerators of preorders}
\author{Yueming Zhong\\[0.5em]
\small Key Laboratory of Engineering Modeling and\\
\small Statistical Computation of Hainan Province\\
\small School of Mathematics and Statistics, Hainan University\\
\small Haikou 570228, P.R. China\\[0.3em]
\small\href{mailto:zhongyueming107@gmail.com}{\texttt{zhongyueming107@gmail.com}}}
\date{}
\hypersetup{pdfauthor={Yueming Zhong},
  pdftitle={Gamma-positivity of support enumerators of preorders}}

\begin{document}
\maketitle

\begin{abstract}
The support enumerator of a finite preorder counts lattice points of its
polytope by support size. We prove Athanasiadis and Chapoton's conjecture
that this polynomial is $\gamma$-positive. On an $n$-element preorder,
the $j$th $\gamma$-coefficient counts points of the preorder polytope
with entries in $\{0,1,2\}$, coordinate sum $n$, and exactly $j$
entries equal to $2$.
Oh's correspondence for transversal matroids and Hall's theorem give
a formula that records the support itself. For arbors, a separate
proof gives a refinement by coordinate sum.
\end{abstract}

\section{Introduction}

Let $\tau$ be a preorder on a finite set $E$ of size $n$. Write
$\preceq_\tau$ for its reflexive, transitive relation.
A subset $I\subseteq E$ is an order ideal if
$r\preceq_\tau s$ and $s\in I$ imply $r\in I$. For a vector $x$ on a set $A$,
write $x(U)=\sum_{i\in U}x_i$ for $U\subseteq A$ and $|x|=x(A)$.
The preorder polytope is
\[
 \Q_\tau=\{x\in\mathbb R_{\ge 0}^{E}:x(I)\le |I|
                       \text{ for every order ideal }I\}.
\]
The ideal $E$ gives the bound $|x|\le n$. For a nonnegative vector $x$,
let $\supp(x)$ be the set of indices $i$ with $x_i>0$. Put
$\bt=(t_i)_{i\in E}$ and $\bt^S=\prod_{i\in S}t_i$ for $S\subseteq E$.
The support enumerators are
\[
 H_\tau(\bt)=\sum_{x\in\Q_\tau\cap\Z^{E}}\bt^{\supp(x)},
 \qquad H_\tau(t)=H_\tau(\bt)\big|_{t_i=t\ (i\in E)} .
\]
The univariate enumerator is denoted by $h(\tau,t)$ in \cite{AC}.

Chapoton introduced the arbor version of the support enumerator \cite{Cha}.
For preorders, Athanasiadis and Chapoton conjectured that the enumerator
is $\gamma$-positive \cite[Conjecture~5.2]{AC}. Cases already established
include linear arbors \cite[Proposition~2.2]{Ath}, octopus arbors
\cite[Theorem~1.3]{AXY}, and lopsided octopuses
\cite[Proposition~3.4]{Menon}. Zig-zag posets and ordinal sums of two
antichains are treated in \cite[Propositions~9.1 and~9.2]{AC}.
Two-row Ferrers posets and recursively defined combs were treated
by Zhong \cite[Theorems~1.2 and~1.3]{Zhong}.

For $a\in\{0,1,2\}^{E}$, put $n_k(a)=|\{i\in E:a_i=k\}|$ and write
\[
 \mathcal A_\tau=\{a\in\Q_\tau\cap\{0,1,2\}^{E}:|a|=n\},
 \qquad
 W_a(\bt)=\prod_{i:a_i=2}t_i\prod_{i:a_i=1}(1+t_i).
\]

\begin{theorem}\label{thm:main}
For every finite preorder $\tau$,
\begin{equation}\label{eq:multi}
 H_\tau(\bt)=\sum_{a\in\mathcal A_\tau}W_a(\bt).
\end{equation}
In particular,
\begin{equation}\label{eq:gamma}
 H_\tau(t)=\sum_{j=0}^{\lfloor n/2\rfloor}
             \gamma_{\tau,j}t^j(1+t)^{n-2j},
 \qquad
 \gamma_{\tau,j}=\bigl|\{a\in\mathcal A_\tau:n_2(a)=j\}\bigr|.
\end{equation}
\end{theorem}

Write $\Gamma_\tau(q)=\sum_{j=0}^{\lfloor n/2\rfloor}\gamma_{\tau,j}q^j$
and set $\gamma_{\tau,j}=0$ for integer indices outside this range.
With zero coefficients outside its degree range, each shifted binomial
polynomial in \eqref{eq:gamma} has coefficients symmetric about degree
$n/2$ and weakly increasing toward that center.
The same holds for their nonnegative sum, proving parts (a) and (b) of
\cite[Conjecture~5.1]{AC}.

The fixed-support count in Proposition~\ref{prop:fixed} uses Oh's
correspondence \cite[Lemma~22 and Proposition~26]{Oh}. It also occurs
in the proof of Ohsugi and Tsuchiya's formula for interior polynomials
\cite[Proposition~3.4]{OT}. In the linear case, Athanasiadis groups pairs
of sets by their differences \cite[Lemma~2.1 and Proposition~2.2]{Ath}.
For general preorders, one Hall criterion on order ideals handles both
the lattice-point assignment and the matching of endpoint sets.
The latter depends only on the two set differences.
Section~\ref{sec:arbors} gives a second proof for arbors and refines
the formula by coordinate sum.

\section{The matching proof}\label{sec:matching}

For $U\subseteq E$, let
\[
 \downarrow U=\{r\in E:r\preceq_\tau s\text{ for some }s\in U\}.
\]
This is an order ideal containing $U$. A pair $(S,B)$ of subsets
of $E$ is \emph{matchable} if a bijection from $B$ to $S$ matches
each $r\in B$ to some $s$ with $r\preceq_\tau s$. Resources and
targets are taken in separate copies of $E$.

We use the following form of Hall's criterion. Given
$x\in\Z_{\ge0}^{E}$ and $B\subseteq E$, assign $x_s$ copies of each
$s\in E$ to distinct resources $r\in B$ with $r\preceq_\tau s$.
Such an assignment exists exactly when
\begin{equation}\label{eq:hall}
 x(I)\le |B\cap I|\qquad(I\text{ an order ideal}).
\end{equation}
Indeed, copies with labels in $I$ can only use resources in $B\cap I$.
Conversely, a set of copies with label set $U$ has at most $x(U)$
members and has neighborhood $B\cap\downarrow U$. Since
$x(U)\le x(\downarrow U)\le |B\cap\downarrow U|$, Hall's theorem
gives the assignment.

\begin{proposition}[Fixed-support count]\label{prop:fixed}
For $S\subseteq E$, the number of lattice points of $\Q_\tau$ with
support $S$ equals the number of sets $B\subseteq E$ for which
$(S,B)$ is matchable.
\end{proposition}

\begin{proof}
For $S=\varnothing$, both counts are one. For $S=E\ne\varnothing$,
the bound $|x|\le n$ leaves only the all-ones vector. The only
possible source set is $B=E$, which is matchable by reflexivity.

Let $\varnothing\ne S\subsetneq E$. The family
$(\downarrow\{s\})_{s\in S}$ presents a transversal matroid on $E$.
Each slot $s$ can use its own label, so the rank is $|S|$ and the
bases are the sets $B$ matchable to $S$. To apply Oh's correspondence,
we put unused resources in a single auxiliary coordinate $\ast$.
The polytope for this presentation is
\[
 \widetilde{\mathcal P}_S
 =\sum_{r\in E}\conv\bigl(\{e_\ast\}\cup
                    \{e_s:s\in S,\ r\preceq_\tau s\}\bigr)
 \subseteq\mathbb R^{S\cup\{\ast\}},
\]
where the $e_i$ are the coordinate vectors.
In the presenting graph, $\ast$ is adjacent to every resource,
and each slot $s$ is adjacent to its own resource. Thus this graph
is connected.

For $x\in\Q_\tau\cap\Z^E$ of support $S$, use \eqref{eq:hall} with
$B=E$ to assign the copies. In summand $r$, choose $e_s$ if $r$ is
assigned to a copy of $s$, and choose $e_\ast$ if $r$ is unused.
This realizes the lattice point $(x|_S,n-|x|)$ of
$\widetilde{\mathcal P}_S$. Conversely, let $y$ be an integral
point of this sum with $y_s\ge1$ for $s\in S$. Extend $y|_S$ by
zero to obtain $x$. For an ideal $I$, only summands with $r\in I$
can contribute to $x(I)$, each by at most one. Hence $x\in\Q_\tau$,
with support $S$. These constructions are inverse because every
point of the sum has coordinate sum $n$.

By Oh's Lemma~22 \cite{Oh}, the lattice points of this sum with
$y_s\ge1$ for $s\in S$ are its base points. Oh's Proposition~26
\cite{Oh} puts them in bijection with the matroid bases.
\end{proof}

\begin{proof}[Proof of Theorem~\ref{thm:main}]
Proposition~\ref{prop:fixed} gives
\begin{equation}\label{eq:endpoints}
 H_\tau(\bt)=\sum_{(S,B)\ \text{matchable}}\bt^S .
\end{equation}
We group the pairs by $a=\ind_E+\ind_S-\ind_B$, where $\ind_A$
is the indicator vector of $A$ on $E$.
For $|S|=|B|$, applying \eqref{eq:hall} to $x=\ind_S$ shows
that $(S,B)$ is matchable exactly when
\begin{equation}\label{eq:ideals}
 |S\cap I|\le |B\cap I|\qquad(I\text{ an order ideal}).
\end{equation}
The entries of $a$ are in $\{0,1,2\}$ and $|a|=n$.
The identity $a(I)=|I|+|S\cap I|-|B\cap I|$ identifies
\eqref{eq:ideals} with $a\in\mathcal A_\tau$.

Fix $a\in\mathcal A_\tau$ and set $P=\{i:a_i=2\}$ and
$N=\{i:a_i=0\}$. Since $|a|=n$, we have $|P|=|N|$.
The pairs giving $a$ are precisely $(P\cup C,N\cup C)$ with
$C\subseteq E\setminus(P\cup N)$. All satisfy \eqref{eq:ideals}.
Here $C=S\cap B$ is determined by the pair. Each pair contributes
once, independently of how it is matched.
The total weight for $a$ is therefore
\begin{equation}\label{eq:class}
 \sum_{C\subseteq E\setminus(P\cup N)}\bt^{P\cup C}
 =\bt^P\prod_{i\in E\setminus(P\cup N)}(1+t_i)
 =W_a(\bt).
\end{equation}
In particular, $C=\varnothing$ gives a bijection between
$\mathcal A_\tau$ and disjoint matchable pairs $(P,N)$.
Summing \eqref{eq:class} proves \eqref{eq:multi}. For $|P|=j$,
the product has $n-2j$ factors. Setting all $t_i=t$ gives
\eqref{eq:gamma}.
\end{proof}

Let $G_\tau$ have resources on the left, targets on the right, and
edges $(r,s)$ for $r\preceq_\tau s$. Write $G_\tau^*$ for the graph with
its two sides exchanged. For a bipartite graph $G$, the augmentation
$\widehat G$ adds a vertex on each side adjacent to every vertex on
the opposite side, including the other new vertex. In the notation
of Dai et al.\ \cite[Theorem~1.2 and Section~4.3]{DHLTW},
\[
 H_\tau(t)=F_{G_\tau^*}(t)=h^*(\mathcal R_{\widehat{G_\tau^*}};t).
\]
Ohsugi and Tsuchiya \cite[Propositions~3.3 and~3.4]{OT} identify
this $h^*$-polynomial with the interior polynomial of
$\widehat{G_\tau^*}$, whose coefficients count matchable endpoint sets
in $G_\tau^*$. Exchanging the two sides preserves these counts, so
the univariate specialization of \eqref{eq:endpoints} also follows
from these results. Dai et al.\ proved invariance under preorder
duality \cite[Theorem~1.3 and Section~4.3]{DHLTW}.

\begin{example}\label{ex:diamond}
Consider the four-element diamond, with a minimum, a maximum, and two
incomparable elements between them. Its five strict comparisons give
$\gamma_{\tau,1}=5$. A target pair must avoid the minimum. The two
middle elements cannot both be targets, as they would need the same
source. If the maximum and one middle element are the targets,
match the minimum to that middle element and the other middle
element to the maximum.
Hence $\gamma_{\tau,2}=2$ and
\[
 H_\tau(t)=(1+t)^4+5t(1+t)^2+2t^2
          =1+9t+18t^2+9t^3+t^4.
\]
\end{example}

\begin{example}[Patel's counterexample to real-rootedness]\label{ex:nonreal}
Patel's eight-element poset \cite{Patel} has lower elements
$a_1,a_2,a_3,d$ and upper elements $b_0,b_1,b_2,b_3$.
Its strict comparisons are $a_i<b_0$, $a_i<b_i$, and $d<b_i$
for $1\le i\le3$. A disjoint matchable pair uses lower elements as
sources and upper elements as targets, so there are nine pairs of
size one.

For targets $\{b_i,b_j\}$ with $i,j>0$, any two sources from
$\{a_i,a_j,d\}$ work. For $\{b_0,b_i\}$, the only excluded source
pair consists of the two $a$'s other than $a_i$. This gives
$3\cdot3+3\cdot5=24$ pairs of size two.
In the full source set, every target has at least two neighbors,
every pair of targets has at least three, and every triple has all four.
Deleting any one source leaves Hall's inequalities for each
target triple, so all sixteen pairs of size three work.
The full target set has all four sources as neighbors, so Hall's
theorem also gives the single pair of size four. Thus
\[
 f(q)=\Gamma_\tau(q)=1+9q+24q^2+16q^3+q^4.
\]
Patel's certificate \cite{Patel} follows from the logarithmic
derivative. If all roots $\rho$ were real, then for real $q$
with $f(q)\ne0$ we would have
\[
 f'(q)^2-f(q)f''(q)
 =f(q)^2\sum_\rho\frac{1}{(q-\rho)^2}\ge0,
\]
where roots are counted with multiplicity. At $q=-1/4$,
$(f,f',f'')=(1/256,-1/16,99/4)$, giving $-95/1024$ on the
left. Hence $f$ has a nonreal root $q_0$.
Since $q_0\ne0$, the equation $t=q_0(1+t)^2$ is quadratic
and has a complex solution. It cannot have $t=-1$ as a solution.
Every real $t\ne-1$ gives a real value of $t/(1+t)^2$, so the
solution is nonreal. Equation \eqref{eq:gamma} gives
$H_\tau(t)=(1+t)^8f(q_0)=0$. This is Patel's counterexample to
\cite[Conjecture~5.3]{AC}.
\end{example}

\section{A height refinement for arbors}\label{sec:arbors}

An arbor $T$ is a rooted tree whose vertices carry nonempty blocks
$B_v$ that partition $E$. Let $D_v$ be the union of the blocks in the
subtree at $v$. Its associated preorder $\tau_T$ has
$r\preceq_{\tau_T}s$ when the vertex carrying $r$ lies in the subtree
at the vertex carrying $s$, including equality of vertices.
Every ideal contains the entire block of each of its labels.
Its maximal blocks in the order induced by $\tau_T$ root pairwise
disjoint subtrees whose label sets partition the ideal. Each $D_v$
is itself an ideal, so the inequalities $x\ge0$ and
$x(D_v)\le |D_v|$ for all $v$ define $\Q_{\tau_T}$.
We use the subscript $T$ in place of $\tau_T$ for this polytope
and its enumerators and coefficients.

\begin{theorem}\label{thm:height}
For an arbor $T$ on a set $E$ of size $n$,
\begin{equation}\label{eq:height}
 \sum_{x\in\Q_T\cap\Z^E}\bt^{\supp(x)}z^{|x|}
 =
 \sum_{a\in\Q_T\cap\{0,1,2\}^{E}}
    W_a(\bt)\,z^{|a|}(1-z)^{n-|a|}.
\end{equation}
\end{theorem}

\begin{proof}
We prove the identity by induction on the number of vertices.
At the root, a three-term decomposition replaces the sum over
nonnegative coordinates by a sum over ternary coordinates.

For a subtree $U$ with label set $E_U$, put $n_U=|E_U|$ and
\[
 \Phi_U(z)=\sum_{x\in\Q_U\cap\Z^{E_U}}\bt^{\supp(x)}z^{|x|}.
\]
Expand all rational factors at $z=0$ as formal power series over
$\Z[t_i:i\in E]$. Write $[\,\cdot\,]_{\le d}$ for truncation
through degree $d$ in $z$.
Let $B$ be the root block of $T$, with child subtrees
$T_1,\ldots,T_m$. Each root coordinate has generating series
$1+t_i z/(1-z)$, while the child factors enforce the proper
subtree bounds. The remaining bound is $|x|\le n$, so
\begin{equation}\label{eq:tree-rec}
 \Phi_T(z)=
 \left[\prod_{i\in B}\left(1+\frac{t_i z}{1-z}\right)
       \prod_{\ell=1}^m\Phi_{T_\ell}(z)\right]_{\le n}.
\end{equation}
An empty child product is $1$. Setting all $t_i=t$ gives
Chapoton's $K_T(t,z)$ and its root-block recursion
\cite[arXiv version, Section~2.5, equation~(16) and Proposition~2.8]{Cha}.

The root factors admit the decomposition
\[
 1+\frac{t_i z}{1-z}
   =(1-z)+(1+t_i)z+\frac{t_i z^2}{1-z}.
\]
The three terms correspond to $a_i=0,1,2$ and have the form
$w_i(a_i)z^{a_i}(1-z)^{1-a_i}$ with weights
$w_i(0)=1$, $w_i(1)=1+t_i$, and $w_i(2)=t_i$.
Use \eqref{eq:height} inductively in the child factors.
Choosing one term per root coordinate and one ternary vector per
child gives a ternary vector $a$ satisfying the proper subtree
bounds. Since $n=|B|+\sum_\ell n_{T_\ell}$, its term is
$W_a(\bt)z^{|a|}(1-z)^{n-|a|}$.

If $|a|\le n$, this term is a polynomial of degree at most $n$
and survives truncation in full. If $|a|>n$, its Taylor expansion
starts in degree $|a|$ and contributes nothing. Thus the
truncation keeps precisely the vectors satisfying the root bound,
giving \eqref{eq:height}. With no child factors the same argument
proves the leaf case.
\end{proof}

At $z=1$, the terms with $|a|<n$ vanish. The terms with $|a|=n$
give \eqref{eq:multi} for arbors.

\begin{corollary}\label{cor:rec}
For an arbor $T$ with label set $E$ of size $n$, put
\[
 C_T(u,q)=\sum_{a\in\Q_T\cap\{0,1,2\}^{E}}
               u^{n-|a|}q^{n_2(a)}.
\]
If the root block has size $r$ and the child subtrees are
$T_1,\ldots,T_m$, then
\begin{equation}\label{eq:dp}
 C_T(u,q)=\Pi_{\ge0}\left[(u+1+q/u)^r
                       \prod_{\ell=1}^{m}C_{T_\ell}(u,q)\right],
 \qquad \Gamma_T(q)=[u^0]C_T(u,q),
\end{equation}
where $\Pi_{\ge0}$ retains the nonnegative powers of $u$.
\end{corollary}

\begin{proof}
Choose a ternary vector satisfying the bounds in each child subtree
and a value in $\{0,1,2\}$ for each root coordinate.
The weight of a root value $0$, $1$, or $2$ is $u$, $1$, or $q/u$,
respectively. The product before truncation in \eqref{eq:dp}
records exactly these choices, with weight $u^{n-|a|}q^{n_2(a)}$
for the resulting vector $a$. The condition $n-|a|\ge0$ is the
remaining root bound. Taking the constant term in $u$ selects
the vectors with $|a|=n$, which give $\Gamma_T(q)$ by \eqref{eq:gamma}.
\end{proof}

\begin{proposition}\label{prop:degree}
For an arbor $T$, let $J_T$ be the graph on $E$ in which two distinct
labels are adjacent if their blocks coincide or one block is a
descendant of the other. Write $\nu(J_T)$ for its maximum matching
size. For every integer $j$,
\[
 \deg\Gamma_T=\nu(J_T),\qquad
 \gamma_{T,j}>0\quad\Longleftrightarrow\quad 0\le j\le\nu(J_T).
\]
The first two coefficients are
\[
 \gamma_{T,0}=1,\qquad
 \gamma_{T,1}=\sum_v |B_v|(|D_v|-1).
\]
\end{proposition}

\begin{proof}
By the proof of Theorem~\ref{thm:main}, $\gamma_{T,j}$ counts
disjoint matchable pairs of size $j$. Any matching realizing such
a pair has $j$ edges with distinct endpoints in $J_T$. Conversely,
orient each edge of a matching in $J_T$ toward the ancestral block,
choosing either direction for an edge within a block. The source
and target sets are disjoint and form a matchable pair.
Thus $\gamma_{T,j}>0$ precisely when $J_T$ has a matching of size
$j$. A maximum matching has submatchings of all sizes from $0$ to
$\nu(J_T)$, which proves the degree and the asserted range.

The empty pair contributes $\gamma_{T,0}=1$. For a target in $B_v$,
there are $|D_v|-1$ choices of a distinct source. Summing over the
targets gives $\gamma_{T,1}$.
\end{proof}

\begin{samepage}
\begin{example}
Take a root block of size two and two isomorphic child arbors on
disjoint label sets. Each child $U$ has a block of size two above
two singleton leaves. A singleton leaf contributes $1+u$, so
\eqref{eq:dp} gives
\[
 C_U(u,q)=\Pi_{\ge0}\bigl[(u+1+q/u)^2(1+u)^2\bigr]
         =(1+u)^4+2q(u^2+3u+3)+q^2.
\]
At the root, \eqref{eq:dp} gives
$\Gamma_T(q)=[u^0](u+1+q/u)^2C_U(u,q)^2$.
The root factor has lowest $u$-degree $-2$, so only the
coefficients $c_i=[u^i]C_U(u,q)$ for $i=0,1,2$ are needed.
Thus
\[
 \begin{aligned}
 \Gamma_T(q)
   &=(1+2q)c_0^2+4qc_0c_1+q^2(c_1^2+2c_0c_2)\\
   &=1+30q+210q^2+372q^3+121q^4+6q^5.
 \end{aligned}
\]
The coefficient of $q$ also follows from
Proposition~\ref{prop:degree}, since the three blocks of size two
contribute $2(10-1)+2(4-1)+2(4-1)$.
\end{example}
\end{samepage}

\end{document}